# On Dual-Finite Volume Methods for Extended Porous Medium Equations

Hidekazu Yoshioka[1]

[1]Faculty of Life and Environmental Science, Shimane University, Shimane, Japan.
E-mail: yoshih@life.shimane-u.ac.jp

## 1. Introduction

Nonlinear conservation laws are mathematically natural models describing transport phenomena of a variety of physical quantities, such as mass, momentum, and energy in continuous medium. Numerically approximating solutions to nonlinear conservation laws requires using a conservative and stable numerical method, so that physically consistent computational results are obtained, such that positivity and/or monotonicity of certain physical quantities are realized. One of the most widely spread numerical methods for the conservation laws is the finite volume method (FVM), which is based on mathematical and numerical local conservation principles. The most important part of the FVM is evaluating numerical fluxes on cell interfaces, which crucially determine accuracy and stability of numerical solutions. For solving parabolic conservation laws with dominant advection in particular, numerical fluxes should guarantee unconditional stability in space but should possess least numerical diffusion as possible. In general, this issue is more delicate for the problems with degenerate coefficients than for those with non-degenerate counterparts. The degeneration of coefficients would cause numerical instability even if sufficiently fine computational meshes were used. Literatures indicate that the mathematical concept of the isotonicity on numerical fluxes, which was originally stated in Ortega and Rheinboldt (1970) and has later been discussed more in detail in Fuhrmann and Langmach (2001), can help develop computationally stable and physically consistent FVMs.

This article focuses on an FVM for numerically solving a class of porous medium equations (PMEs), which govern moisture dynamics in porous media, such as soils (Szymkiewicz, and Helmig, 2011; Zha et al., 2013) and thin fibrous media (Landeryou et al., 2005; Lockington et al., 2007; Ashari et al., 2010). Similar differential equations are encountered in different physical problems, such as heat conduction problems in plasma (Bertsch, 1982; Wilhelmsson, 1988), combustion of liquids (Lakkis and Ghoniem, 2003; Rida, 2010), evaporation dynamics of volatile liquids (Okrasisk et al., 2001; Rida, 2010), chemotaxis of cells and organisms (Borsche et al., 2014), and astrophysics (Haerns and Van Gorder, 2012). Due to their importance in science and engineering applications, mathematical and numerical analyses on the PMEs have enthusiastically been performed in a vast number of researches so far (Antontsev et al., 2002; Antontsev and Shemarev, 2015; Vázquez, 2007).

Finite volume numerical modeling of the PMEs and related differential equations have focused on computing numerical solutions that are physically consistent. Toro and Hidalgo (2009) proposed a low to very high order FVMs for approximating solutions to 1-D diffusion-reaction equations that generalizes PMEs, focusing on analytically and numerically consistent treatment of the source terms and fluxes. Zamba et al. (2012) proposed a computationally accurate oscillation free FVM based on a weighted essentially non-oscillatory interpolation functions. The Discontinuous Galerkin methods (DGMs) are finite volume analogues of finite element methods (FEMs), which can also achieve mass conservation properties (Yan, 2013; Zhang et al., 2009). Finite Volume Element Methods (FVEMs) are conservative numerical methods combining finite element discretization and finite volume numerical flux evaluation, which have also been effectively used for solving the PMEs in the past researches (Cumming et al., 2011; Albets-Chiko and

Kassions, 2013).

According to Fuhrmann and Langmach (2001), numerical fluxes in the FVMs for nonlinear differential equations such as PMEs should equip with the isotonicity, which is rigorously defined later, for computing physically consistent numerical solutions without using extremely fine computational meshes. Their numerical fluxes for the PMEs without advection terms are isotone and thus guarantee the above-mentioned mathematical properties of numerical solutions. However, the numerical fluxes that they presented would not be optimal for the PMEs with advection terms typically encountered in real applications as demonstrated later in this article.

The main purpose of this article is to propose a computationally stable FVM for PMEs, focusing in particular on its applications to the problems encountered in physically modelling moisture dynamics in non-woven fibrous sheets under evaporative environments. The governing equation in this case is given by a PME having an evaporation term, which has merely been encountered in the literatures. The FVM is developed on the basis of the Dual-Finite Volume Method (DFVM) proposed in Yoshioka and Unami (2013) originally proposed for simulating solute dispersion phenomena in surface water networks modelled as connected graphs consisting of 1-D reaches and 0-D junctions. The DFVM is unconditionally stable in space for linear parabolic equations and is unconditionally stable in both space and time if equipped with an appropriate implicit temporal integration method. Numerical fluxes are evaluated using analytical solutions to local two-point boundary value problems in the DFVM, and this numerical technique is referred to as the fitting technique. Numerical methods using the fitting technique have successfully been used for solving differential equations encountered in a variety of problems. Such examples include the Black-Scholes equation governing option pricing (Chernogorova and Valkov, 2011; Wang et al., 2014), the Hamilton-Jacobi-Bellman equation appearing in optimization problems (Richardson and Wang, 2006), and the reactive solute transport equation in incompressible fluids (Fuhrmann et al., 2011; Fiebach et al., 2014). However, to the author's knowledge, the fitting technique has not been applied to the PMEs with the evaporation term, which is a main motivation of this article. In addition, it is also scientifically and practically interesting to analyze solutions to the PMEs on connected graphs that is approached in this article from a numerical point of view, which is encountered in the problems of moisture dynamics in fibrous sheet networks.

This article shows that the unconditional stability of the DFVM, which is at least valid for linear problems, is not true for generic nonlinear differential equations including the PMEs unless the coefficient appearing in the numerical fluxes are appropriately evaluated. Unconditional stability of the DFVM for the PMEs can be achieved using the concept of the isotonicity in this article. This article shows that the fully-upwind type coefficient proposed in Fuhrmann and Langmach (2001), which has widely been used in practical problems because of its convenience and stability (Belfort et al., 2013; Lai and Ogden, 2015), is not truly isotonic for our PMEs. This article provides a theoretically truly isotone numerical fluxes specialized for solving the PMEs presented, which is still as simple as the conventional fully-upwind counterpart.

The rest of this article is organized as follows. Section 2 presents the 1-D and 2-D PMEs governing moisture dynamics in thin non-woven fibrous sheets. Their regularized counterparts for well-posing the problems under wider conditions are also presented in this section. Section 3 presents spatial and temporal discretization procedures of a 1-D DFVM with the particular emphasis on evaluating of the numerical fluxes with the isotonicity condition. Section 4 presents a 2-D counterpart of the 1-D DFVM, which is based on analogous spatial and temporal discretization procedures.

## 2. Mathematical Model

The 1-D and 2-D PMEs for mathematically describing moisture dynamics in thin non-woven fibrous sheets subject to evaporation are presented. Their regularized counterparts are also presented in this section.

### 2.1 1-D Porous Medium Equation

The 1-D PME is presented for simulating 1-D moisture dynamics in porous media, such as longitudinal moisture dynamics subject to evaporation in inclined thin non-woven fibrous sheets. It is reasonable to assume that moisture profiles in the sheets are transversely homogenous, leading to a longitudinally 1-D mathematical model. The saturation at each longitudinally 1-D position $x$ of a sheet at the time $t$ is denoted by $\theta = \theta(t,x)$, and is normalized as $0 \leq u = (\theta - \theta_r)(\theta_s - \theta_r)^{-1} \leq 1$ where $\theta_M (=0.9)$ and $\theta_r (=0.0)$ are the maximum and minimum water contents of the sheet. Based on the non-dimensionalization method, the extended PME is presented as (Landeryou et al., 2005; Lockington et al., 2007)

$$\frac{\partial u}{\partial \tilde{t}} = \frac{\partial}{\partial \tilde{x}}\left((m-p)u^{m-p-1}\frac{\partial u}{\partial \tilde{x}}\right) + \sin\alpha\frac{\partial}{\partial \tilde{x}}\left(u^{m-1}u\right) - \tilde{E}_s u^q \tag{1}$$

with the non-dimensional independent variables

$$\tilde{t} = \frac{K_s^2(m-p)}{D_s \theta_s^2}t, \quad \tilde{x} = \frac{K_s(m-p)}{D_s \theta_s}x, \text{ and } \tilde{E}_s = \frac{D_s \theta_s^2}{K_s^2(m-p)}E_s \tag{2}$$

where $D_s(>0)$ is the saturated diffusivity, $K_s(>0)$ is the saturated permeability, $E_s(\geq 0)$ is the evaporation coefficient, $m$, $p(<m)$, and $q$ are the positive parameters determining nonlinearity of Eq.(1), and $-0.5\pi \leq \alpha \leq 0.5\pi$ is the inclination angle of the sheet. A remarkable difference between the extended PME(1) and the conventional counterparts is that the former can simulate the evaporative moisture dynamics that the latter in principle cannot do. In deriving Eq.(1), the pressure head $\psi(u)$ and the permeability $K(u)$ are physically parameterized as

$$\psi(u) = \frac{\theta_s D_s}{pK_s}\left(1 - \frac{1}{u^p}\right) \text{ and } K(u) = K_s u^m \tag{3}$$

respectively, which reduce to the conventional models for $p \to +0$ (Brooks and Corey, 1964; Smiles, 1998; Landeryou et al., 2005; Lockington et al., 2007; Ashari et al., 2010). The angle $\alpha$ may be spatially distributed. There exist two contrasting cases on arrangement of the sheets, which are the horizontal ($\alpha = 0$) and vertical cases ($\alpha = \pm 0.5\pi$). The advection term of Eq.(1) vanishes in the horizontal case because $\sin\alpha = 0$. Hereafter, "~" representing non-dimensional variables are omitted from the variables.

Degeneration of the advection and diffusion terms of Eq.(1) would severe mathematical and computational difficulties due to their low regularities. The degeneration can be effectively mitigated with a regularization technique, which regularizes Eq.(1) as

$$\frac{\partial u}{\partial t} = \frac{\partial}{\partial x}\left((m-p)\left[h_\varepsilon(u)\right]^{m-p-1}\frac{\partial u}{\partial x}\right) + \sin\alpha\frac{\partial}{\partial x}\left(\left[h_\varepsilon(u)\right]^{m-1}u\right) - E_s u^q \tag{4}$$

using the regularization kernel $h_\varepsilon(u) = \sqrt{u^2 + \varepsilon^2} (\geq \varepsilon)$ with a small positive parameter $\varepsilon(=10^{-10})$, which completely rules out the degeneration of the advection and diffusion terms. An application of some of the conventional regularization method to Eq.(1) remains as an another option to be applied to Eq.(1); however, they rely on variable transformations of the solution $u$ to a new variable, which in general cannot numerically realize mass conservation property (Pop et al., 2004).

Eq.(4) can be rewritten in the flux form as

$$\frac{\partial u}{\partial t} + \frac{\partial F}{\partial x} = -E_s u^q \tag{5}$$

with the flux $F = F(u)$ given by

$$F = -(m-p)\left[h_\varepsilon(u)\right]^{m-p-1}\frac{\partial u}{\partial x} - \sin\alpha\left[h_\varepsilon(u)\right]^{m-1} u. \tag{6}$$

## 2.2 2-D Porous Medium Equation

The 2-D PME, which is a 2-D counterpart of Eq.(1), is presented in this sub-section. This mathematical model is suited to simulating essentially 2-D moisture dynamics in inclined fibrous sheets. The inclination angle is again denoted by $\alpha$, which typically represents horizontal ($\alpha = 0$) and vertical arrangements ($\alpha = \pm 0.5\pi$) of the sheet. The 2-D $x_1 - x_2$ Cartesian coordinates is taken in the domain $\Omega$. It is assumed that the sheet is horizontal on the $x_2$-direction. Under this condition, the 2-D PME governing the spatio-temporal evolution of the saturation $u = u(t, x_1, x_2)$ is given by (Landeryou et al., 2005; Lockington et al., 2007)

$$\frac{\partial u}{\partial t} = \frac{\partial}{\partial x_1}\left((m-p)u^{m-p-1}\frac{\partial u}{\partial x_1}\right) + \frac{\partial}{\partial x_2}\left((m-p)u^{m-p-1}\frac{\partial u}{\partial x_2}\right) + \sin\alpha\frac{\partial}{\partial x_1}\left(u^{m-1}u\right) - E_s u^q \tag{7}$$

where the known functions are analogously defined as in the 1-D case. Following Eq.(4), The regularized counterpart of Eq.(7) is given by

$$\frac{\partial u}{\partial t} = \frac{\partial}{\partial x_1}\left((m-p)\left[h_\varepsilon(u)\right]^{m-p-1}\frac{\partial u}{\partial x_1}\right) + \frac{\partial}{\partial x_2}\left((m-p)\left[h_\varepsilon(u)\right]^{m-p-1}\frac{\partial u}{\partial x_2}\right)$$
$$+ \sin\alpha\frac{\partial}{\partial x_1}\left(\left[h_\varepsilon(u)\right]^{m-1} u\right) - E_s u^q \tag{8}$$

Eq.(8) can be rewritten in the flux form as

$$\frac{\partial u}{\partial t} + \frac{\partial F_1}{\partial x_1} + \frac{\partial F_2}{\partial x_2} = -E_s u^q \tag{9}$$

with the fluxes $(F_1, F_2)$ given by

$$\begin{pmatrix} F_1 \\ F_2 \end{pmatrix} = \begin{pmatrix} -(m-p)\left[h_\varepsilon(u)\right]^{m-p-1}\frac{\partial u}{\partial x_1} - \sin\alpha\left[h_\varepsilon(u)\right]^{m-1} u \\ -(m-p)\left[h_\varepsilon(u)\right]^{m-p-1}\frac{\partial u}{\partial x_2} \end{pmatrix}. \tag{10}$$

## 3.  1-D Dual-Finite Volume Method

A numerical method for approximating solutions to Eq.(4), which is referred to as the 1-D Dual-Finite Volume Method (DFVM), is explained in this section. The 1-D DFVM was originally developed for solving the Kolmogorov's forward equations, which are linear and conservative advection-diffusion-decay equations (ADDEs), on connected graphs (Yoshioka and Unami, 2013). The 1-D DFVM uses the primal and dual computational meshes and evaluates the numerical fluxes based on analytical solutions to local two-point boundary value problems. Mathematical analysis revealed that the 1-D DFVM is unconditionally stable with an appropriate implicit temporal integration algorithm.

### 3.1 Computational Mesh

Firstly, the domain $\Omega$, which is assumed to be given by a connected graph domain in this sub-section, is divided into a regular mesh that consists of regular cells bounded by two nodes, so that any vertex in the domain, such as a junction, falls on one of the nodes. The regular cells and the nodes are indexed with the natural numbers. The total numbers of regular cells and nodes are denoted by $N_c$ and $N_n$, respectively. The $i$ th node in the domain $\Omega$ is denoted by $P_i$ with its $x$ abscissa $x_i$. The $k$ th regular cell is denoted by $\Omega_k$. The length of $\Omega_k$ is represented by $l_k$. The two nodes bounding both sides of $\Omega_k$ are denoted by the $\varphi(k,1)$ th node $P_{\varphi(k,1)}$ and the $\varphi(k,2)$ th node $P_{\varphi(k,2)}$. The $x$ abscissa in the cell $\Omega_k$ is directed from the node $P_{\varphi(k,1)}$ to the node $P_{\varphi(k,2)}$. The number of regular cells connected to the node $P_i$ is denoted by $\nu(i)$. The $j$ th regular cell connected to the node $P_i$ is referred to as the $\kappa(i,j)$ th regular cell $\Omega_{\kappa(i,j)}$. There exist two nodes that bound the cell $\Omega_{\kappa(i,j)}$; one is the $i$ th node $P_i$ and the other is referred to as the $\mu(i,j)$ th node $P_{\mu(i,j)}$. In the regular cell $\Omega_{\kappa(i,j)}$, the direction of the abscissa is identified with the parameter $\sigma_{i,j}$, which is equal to $1$ when $x$ is directed to the node $P_i$, and is otherwise equal to $-1$. A dual mesh is generated from the regular mesh. Following the multi-dimensional analogue (Mishev, 1998), the $i$ th dual cell $S_i$ is associated with the node $P_i$, which is defined as

$$S_i = \left\{ x \mid |x_i - x| < |x_{\mu(i,j)} - x| \text{ for } 1 \leq j \leq \nu(i) \right\}. \tag{11}$$

The dual mesh consists of $N_n$ dual cells. The cell interface between the dual cells $S_i$ and $S_{\mu(i,j)}$ is denoted by $\Gamma_{i,j}$. The solution $u$ is attributed to the dual cells, and the flux $F$ is attributed to the regular cells.

### 3.2 Operator-Splitting Algorithm

In the present 1-D DFVM, Eq.(4) is regarded as an advection-diffusion-decay equation having solution-dependent coefficients, which are evaluated fully explicitly at each time step. An operator-splitting algorithm analogous to that of Li et al. (2010) specialized for solving Eq.(4) is developed and incorporated into the 1-D DFVM, so that computational instability due to nonlinearity of the evaporation term is completely avoided. The time increment at each time step as denoted by $\Delta t$ and the differential operators defining the advection and diffusion terms of Eq. (4) as $P_{AD}$ and that for the evaporation term as $P_E$. The 1-D DFVM solves Eq.(4) at each time step as

$$u^{(k+1)} = \exp(0.5\Delta t P_E)\exp(\Delta t P_{AD})\exp(0.5\Delta t P_E)u^{(k)}, \tag{12}$$

where $u^{(k)}$ is the solution at the $k$ th time step. Although the 1-D DFVM is unconditionally stable for linear problems, preliminary computation showed that it is not the case for the PMEs with advection terms if $\Delta t$ is sufficiently large. The reason of this issue would be the loss of ellipticity of the resulting discretized system with large $\Delta t$ as implied in Pop et al. (2004). This drawback can be overcome if an implicit treatment of the coefficients is used, which will be addressed in future researches. Discretization procedure of the two sub-steps at each time step is explained

### 3.3 Evaporation Sub-step

At the evaporation sub-step, the equation to be solved at each node is given by

$$\frac{\partial u}{\partial t} = -E_s u^q, \tag{13}$$

which can be formally regarded as an ordinary differential equation (ODE). At the $i$ th node, starting from an initial condition $U_i$ at the time $t_0$, the analytical nodal solution to Eq.(13) at the time $t > t_0$ is analytically expressed as

$$u_i(t) = \begin{cases} \left(\max\{0, \phi_i(U_i, t, t_0)\}\right)^{\frac{1}{1-q}} & (0 < q < 1) \\ U_i \exp(-E_s t) & (q = 1) \\ \phi_i(U_i, t, t_0)^{-\frac{1}{q-1}} & (q > 1) \end{cases} \quad (14)$$

with

$$\phi_i(u, t, t_0) = u^{1-q} - (1-q) E_s (t - t_0) \quad (15)$$

for $u \in \mathbb{R}^+$ without any numerical approximations. By Eq.(14), the first evaporation sub-step at the $k$ th time step for the $i$ th node is mathematically expressed as

$$u_i^* = \exp_d(0.5\Delta t P_E) u^{(k)} \quad (16)$$

with

$$\exp_d(0.5\Delta t P_E) u_i^{(k)} = \begin{cases} \left(\max\{0, \phi_i(u_i^{(k)}, t + \Delta t, t)\}\right)^{\frac{1}{1-q}} & (0 < q < 1) \\ u_i^{(k)} \exp(-E_s \Delta t) & (q = 1) \\ \phi_i(u_i^{(k)}, t + \Delta t, t)^{-\frac{1}{q-1}} & (q > 1) \end{cases} \quad (17)$$

where $u_i^*$ represents the updated nodal value, which is used as the initial condition of the advection-diffusion sub-step. The updated nodal value at the $i$ th node just after the advection-diffusion sub-step is hereafter denoted by $u_i^{**}$. The second evaporation sub-step is mathematically expressed as

$$u^{(k+1)} = \exp_d(0.5\Delta t P_E) u_i^{**} \quad (18)$$

with

$$\exp_d(0.5\Delta t P_E) u_i^{**} = \begin{cases} \left(\max\{0, \phi_i(u_i^{**}, t + 0.5\Delta t, t)\}\right)^{\frac{1}{1-q}} & (0 < q < 1) \\ u_i^{**} \exp(-0.5 E_s \Delta t) & (q = 1) \\ \phi_i(u_i^{**}, t + 0.5\Delta t, t)^{-\frac{1}{q-1}} & (q > 1) \end{cases}. \quad (19)$$

It should be noted that there is no restriction on the time increment $\Delta t$ at the evaporation sub-steps because the temporal integration is performed exactly based on the knowledge on the analytical solution.

### 3.4 Advection-diffusion Sub-step

At the advection-diffusion sub-step, the equation to be solved is given by

$$\frac{\partial u}{\partial t} = \frac{\partial}{\partial x}\left((m - p)[h_\varepsilon(u)]^{m-p-1}\frac{\partial u}{\partial x}\right) + \sin\alpha \frac{\partial}{\partial x}\left([h_\varepsilon(u)]^{m-1} u\right), \quad (20)$$

which can be formally regarded as a nonlinear advection-diffusion equation. A finite volume discretization is applied to Eq.(20) as presented in the following sub-sections.

#### 3.4.1 Finite Volume Formulation

Integrating Eq.(20) in the dual cell $S_i$ with the application of the Gauss-Green theorem leads to the local integral relationship

$$\frac{\partial}{\partial t}\int_{S_i} u \, ds + \sum_{j=1}^{\nu(i)} \sigma_{i,j} F_{i,j} = 0. \quad (21)$$

where $F_{i,j}$ is the numerical flux to be evaluated on the cell interface $\Gamma_{i,j}$. The integrals in Eq.(21) prescribe an internal boundary condition when the dual cell $S_i$ contains a junction, overcoming the difficulties that the most numerical schemes, such as Szymkiewicz (2008) and Zhang et al. (2010), encounter. The first term in the left hand-side of Eq.(21) is evaluated as

$$\frac{\partial}{\partial t}\int_{S_i} u \, dx \approx |S_i|\frac{du_i}{dt}, \tag{22}$$

where $|S_i|$ is the measure of the dual cell $S_i$. Elementary calculation leads to that the measure $|S_i|$ is given by

$$|S_i| = \frac{1}{2}\sum_{j=1}^{v(i)} l_{\kappa(i,j)}. \tag{23}$$

The numerical flux $F_{i,j}$ on the cell interface $\Gamma_{i,j}$ is evaluated with a fitting technique, which uses the analytical solution $\tilde{u}$ of the linearized two-point boundary value problem

$$\frac{\partial}{\partial x}\left(a_{\kappa(i,j)} u - b_{\kappa(i,j)}\frac{\partial u}{\partial x}\right) = 0 \text{ in } \Omega_{\kappa(i,j)} \tag{24}$$

subject to the nodal boundary conditions

$$u(x_i) = u_i, \quad u(x_{\mu(i,j)}) = u_{\mu(i,j)} \tag{25}$$

where the coefficients $a_{\kappa(i,j)}$ and $b_{\kappa(i,j)}$ in Eq.(24) are given by

$$a_{\kappa(i,j)} = (m-p)\left[h_\varepsilon\left(\overline{u_{\kappa(i,j)}}\right)\right]^{m-p-1} \tag{26}$$

and

$$b_{\kappa(i,j)} = -\sin\alpha\left[h_\varepsilon\left(\overline{u_{\kappa(i,j)}}\right)\right]^{m-1}, \tag{27}$$

respectively where $\overline{u_{\kappa(i,j)}} = \overline{u_{\kappa(i,j)}}(u_i, u_{\mu(i,j)})$ represents the value of $u$ in the regular cell $\Omega_{\kappa(i,j)}$, which turns out to be a critical factor determining computational accuracy and stability of the DFVM for the PME, and it is specified later. The function $\overline{u_{\kappa(i,j)}}$ at least satisfies

$$\overline{u_{\kappa(i,j)}}(u,u) = u \tag{28}$$

for arbitrary $u \in \mathbb{R}$. The analytical solution $u$ to the cell-wise two-point boundary value problem is explicitly obtained as

$$u = \frac{\exp(p_{i,j})u_i - u_{\mu(i,j)}}{\exp(p_{i,j}) - 1} + \frac{u_{\mu(i,j)} - u_i}{\exp(p_{i,j}) - 1}\exp\left(\frac{\sigma_{i,j} b_{\kappa(i,j)}}{a_{\kappa(i,j)}}(x - x_i)\right), \tag{29}$$

which is used for directly evaluating the numerical flux $F_{i,j}$ as

$$\begin{aligned} F_{i,j} &= a_{\kappa(i,j)} u - b_{\kappa(i,j)}\frac{\partial u}{\partial x} \\ &= \frac{\sigma_{i,j} a_{\kappa(i,j)}}{l_{\kappa(i,j)}} p_{i,j}\frac{\exp(p_{i,j})u_i - u_{\mu(i,j)}}{\exp(p_{i,j}) - 1} \\ &= \frac{(m-p)\sigma_{i,j}}{l_{\kappa(i,j)}}\left[h_\varepsilon\left(\overline{u_{\kappa(i,j)}}\right)\right]^{m-p-1} p_{i,j}\frac{\exp(p_{i,j})u_i - u_{\mu(i,j)}}{\exp(p_{i,j}) - 1} \end{aligned} \tag{30}$$

where the local Peclet number $p_{i,j}$, which in the present case depends on the solution $u$, is defined as

$$p_{i,j} = \frac{\sigma_{i,j} b_{\kappa(i,j)} l_{\kappa(i,j)}}{a_{\kappa(i,j)}} = -\frac{\sigma_{i,j} l_{\kappa(i,j)} \sin\alpha}{m-p} h_\varepsilon\left(\overline{u_{\kappa(i,j)}}\right)^p. \tag{31}$$

The numerical flux $F_{i,j}$ is fully specified if the value of $\overline{u_{\kappa(i,j)}}$ is evaluated. It should be noted that evaluation of the numerical flux using the fitting technique ensures at least first-order spatial convergence of the scheme for linear problems (Roos, 1994; Wang, 2004). A remarkable point of the present DFVM is that the Peclet number $p_{i,j}$ is formally independent of the solution $u$ when $p=0$. In fact, by Eq.(31), $p_{i,j}$ in this case reduces to

$$p_{i,j} = -\frac{\sigma_{i,j} l_{\kappa(i,j)} \sin\alpha}{m-p}. \tag{32}$$

This remarkable property help propose a truly isotonic numerical flux for the special but still realistic case with $p=0$ as presented in the next section.

### 3.4.2  Numerical Flux

The value of $\overline{u_{\kappa(i,j)}}$ in the element $\Omega_{\kappa(i,j)}$ is specified for completely determining the numerical flux $F_{i,j}$. The three evaluation methods of $\overline{u_{\kappa(i,j)}}$ are presented in this article, which are the central (CE) method, fully-upwind (FU) method, and the isotonic (IS) method. Both the CE and FU have been used in the conventional researches, and the IS method is a new numerical method that is proposed in this article. For the sake of simplicity of calculation, the approximation $\varepsilon=0$ is used in this section because it is assumed to be given as a sufficiently small positive value such that it does not significantly affect computational accuracy of numerical solutions. The CE method and the FU method evaluate $\overline{u_{\kappa(i,j)}}$ as

$$\overline{u_{\kappa(i,j)}} = \frac{1}{2}\left(u_i + u_{\mu(i,j)}\right) \tag{33}$$

and

$$\overline{u_{\kappa(i,j)}} = \max\left(u_i, u_{\mu(i,j)}\right), \tag{34}$$

respectively.

The concept of isotonicity is presented in this sub-section for proposing an isotonic numerical flux suited to solving the PME. The numerical flux is regarded as a bi-variate function of the parameters $\lambda, \mu \geq 0$ as

$$F_{i,j} = F_{i,j}(\lambda,\mu)\big|_{\lambda=u_i,\mu=u_{\mu(i,j)}} = \frac{\sigma_{i,j} p_{i,j}}{l_{\kappa(i,j)}}\left[a_{\kappa(i,j)}(\lambda,\mu)\frac{\exp(p_{i,j})\lambda - \mu}{\exp(p_{i,j})-1}\right]_{\lambda=u_i,\mu=u_{\mu(i,j)}}. \tag{35}$$

In the present case, the numerical flux $F_{i,j}$ is said to be isotone if the conditions

$$\sigma_{i,j}\frac{\partial F_{i,j}}{\partial \lambda} \geq 0 \tag{36}$$

and

$$-\sigma_{i,j}\frac{\partial F_{i,j}}{\partial \mu} \geq 0 \tag{37}$$

are satisfied for arbitrary $\lambda, \mu \geq 0$ (Fuhrmann and Langmach, 2001). The value of $\overline{u_{\kappa(i,j)}}$ should be determined so that both Eqs.(36) and (37) are rigorously satisfied. However, analytically finding such $\overline{u_{\kappa(i,j)}}$ is difficult because of the nonlinearity of the numerical flux $F_{i,j}$. An exceptional case occurs when $p=0$ where the analytical expression of $\overline{u_{\kappa(i,j)}}$ can be found as demonstrated in the next subsection.

### 3.4.3  Isotonicity on the CE Method

The CE method would give unphysical (negative) numerical solutions with sudden transitions of solution profiles because it does not comply with Eqs.(36) and (37) even with the absence of advection term (Fuhrmann and Langmach, 2001).

### 3.4.4 Isotonicity on the FU Method

The FU method can potentially mitigate such computational difficulty but would give numerical solutions with excessive numerical diffusion. The literature indicated that the FU method is isotone for the problems without the advection; however, this is not true for the problems with advection. For $\lambda \geq \mu \geq 0$ with the FU scheme, the quantities $\sigma_{i,j} \frac{\partial F_{i,j}}{\partial \lambda}$ and $-\sigma_{i,j} \frac{\partial F_{i,j}}{\partial \mu}$ are calculated as

$$\begin{aligned}
\sigma_{i,j} \frac{\partial F_{i,j}}{\partial \lambda} &= \frac{\partial}{\partial \lambda}\left[\frac{mp_{i,j}}{l_{\kappa(i,j)}} \lambda^{m-1} \frac{\exp(p_{i,j})\lambda - \mu}{\exp(p_{i,j})-1}\right] \\
&= \frac{m}{l_{\kappa(i,j)}} \frac{p_{i,j}}{\exp(p_{i,j})-1} \frac{\partial}{\partial \lambda}\left[\lambda^{m-1}\left(\exp(p_{i,j})\lambda - \mu\right)\right] \\
&= \frac{m}{l_{\kappa(i,j)}} \frac{p_{i,j}}{\exp(p_{i,j})-1}\left(m\exp(p_{i,j})\lambda^{m-1} - (m-1)\lambda^{m-2}\mu\right) \\
&= \frac{m}{l_{\kappa(i,j)}} \frac{p_{i,j}}{\exp(p_{i,j})-1} \lambda^{m-2}\left(m\exp(p_{i,j})\lambda - (m-1)\mu\right)
\end{aligned} \qquad (38)$$

and

$$\begin{aligned}
-\sigma_{i,j} \frac{\partial F_{i,j}}{\partial \mu} &= \frac{\partial}{\partial \mu}\left[\frac{mp_{i,j}}{l_{\kappa(i,j)}} \lambda^{m-1} \frac{\mu - \exp(p_{i,j})\lambda}{\exp(p_{i,j})-1}\right] \\
&= \frac{m\lambda^{m-1}}{l_{\kappa(i,j)}} \frac{p_{i,j}}{\exp(p_{i,j})-1} \\
&\geq 0
\end{aligned} \qquad (39)$$

respectively. Similarly, for $\mu \geq \lambda \geq 0$ with the FU scheme, the inequalities are calculated as

$$\begin{aligned}
\sigma_{i,j} \frac{\partial F_{i,j}}{\partial \lambda} &= \frac{\partial}{\partial \lambda}\left[\frac{mp_{i,j}}{l_{\kappa(i,j)}} \mu^{m-1} \frac{\exp(p_{i,j})\lambda - \mu}{\exp(p_{i,j})-1}\right] \\
&= \frac{m}{l_{\kappa(i,j)}} \frac{p_{i,j}}{\exp(p_{i,j})-1} \frac{\partial}{\partial \lambda}\left[\mu^{m-1}\left(\exp(p_{i,j})\lambda - \mu\right)\right] \\
&= \frac{m}{l_{\kappa(i,j)}} \frac{p_{i,j}\exp(p_{i,j})}{\exp(p_{i,j})-1} \mu^{m-1} \\
&\geq 0
\end{aligned} \qquad (40)$$

and

$$\begin{aligned}
-\sigma_{i,j}\frac{\partial F_{i,j}}{\partial \mu} &= \frac{\partial}{\partial \mu}\left[\frac{mp_{i,j}}{l_{\kappa(i,j)}}b_{\kappa(i,j)}\frac{\mu-\exp(p_{i,j})\lambda}{\exp(p_{i,j})-1}\right]\\
&= \frac{m}{l_{\kappa(i,j)}}\frac{p_{i,j}}{\exp(p_{i,j})-1}\frac{\partial}{\partial \mu}\left[\mu^{m-1}\left(\mu-\exp(p_{i,j})\lambda\right)\right]\\
&= \frac{m}{l_{\kappa(i,j)}}\frac{p_{i,j}}{\exp(p_{i,j})-1}\left[m\mu^{m-1}-(m-1)\exp(p_{i,j})\mu^{m-2}\lambda\right]\\
&= \frac{m}{l_{\kappa(i,j)}}\frac{p_{i,j}\mu^{m-2}}{\exp(p_{i,j})-1}\left[m\mu-(m-1)\exp(p_{i,j})\lambda\right]
\end{aligned} \qquad (41)$$

respectively. By Eq.(38), the condition (36) is true if

$$\exp(p_{i,j})-\frac{m-1}{m}\geq 0. \qquad (42)$$

Similarly, by Eq.(41), the condition (37) is true if

$$\exp(-p_{i,j})-\frac{m-1}{m}\geq 0. \qquad (43)$$

Both of the inequality is true if

$$\frac{m-1}{m}\leq \exp(p_{i,j})\leq \frac{m}{m-1}. \qquad (44)$$

If $m=1$, the PME is linear and Eq.(44) reduces to

$$0\leq \exp(p_{i,j})\leq +\infty, \qquad (45)$$

which is a trivially-satisfied inequality. The above-presented mathematical analysis demonstrates that the FU scheme is theoretically not isotone when the inequality (44) is violated.

### 3.4.5  Isotonicity on the IS Method

A new evaluation method of the numerical flux, which is the IS method, is presented as a possible alternative to the CE and FU method. The IS method is expressed for $p=0$ as

$$\overline{u_{\kappa(i,j)}} = \begin{cases} u_i & \left(\text{if } \exp(p_{i,j})u_i - u_{\mu(i,j)} \geq 0\right)\\ u_{\mu(i,j)} & (\text{Otherwise}) \end{cases}. \qquad (46)$$

By Eqs.(38) and (46), isotonicity of the IS method can be proven for the case $\exp(p_{i,j})u_i - u_{\mu(i,j)} \geq 0$ as

$$\begin{aligned}
\sigma_{i,j}\frac{\partial F_{i,j}}{\partial \lambda} &= \frac{m}{l_{\kappa(i,j)}}\frac{p_{i,j}}{\exp(p_{i,j})-1}\lambda^{m-2}\left(m\exp(p_{i,j})\lambda - (m-1)\mu\right)\\
&\geq \frac{m}{l_{\kappa(i,j)}}\frac{p_{i,j}}{\exp(p_{i,j})-1}\lambda^{m-2}\left(m\exp(p_{i,j})\lambda - (m-1)\exp(p_{i,j})\lambda\right),\\
&= \frac{m}{l_{\kappa(i,j)}}\frac{p_{i,j}\exp(p_{i,j})}{\exp(p_{i,j})-1}\lambda^{m-1}\\
&\geq 0
\end{aligned} \qquad (47)$$

which verifies Eq.(36). Similarly, Eq.(37) is verified by Eqs.(41) and (46) as

$$-\sigma_{i,j}\frac{\partial F_{i,j}}{\partial \mu} = \frac{m}{l_{\kappa(i,j)}}\frac{p_{i,j}\mu^{m-2}}{\exp(p_{i,j})-1}\left[m\mu-(m-1)\exp(p_{i,j})\lambda\right]$$
$$\geq \frac{m}{l_{\kappa(i,j)}}\frac{p_{i,j}\lambda\mu^{m-2}}{\exp(p_{i,j})-1}\left[m\mu-(m-1)\mu\right]$$
$$= \frac{m}{l_{\kappa(i,j)}}\frac{p_{i,j}\lambda\mu^{m-1}}{\exp(p_{i,j})-1}$$
$$\geq 0$$
(48)

Eqs.(39), (40), (47), and (48) prove that the IS scheme is isotone.

### 3.4.6 Temporal Integration

Assembling Eq.(21) with the above-mentioned spatial discretization procedure for every node yields the system of linear ODEs. The DFV scheme finally yields the system of ordinary differential equations (ODEs) of the form

$$\frac{d\mathbf{u}}{dt} = \Xi\mathbf{u} + \mathbf{d}$$
(49)

where $\mathbf{u}=[u_i]$ is the $N_n$-dimensional solution vector, $\Xi=[\Xi_{i,k}]$ is the $N_n \times N_n$-dimensional matrix arising from the spatial discretization, and $\mathbf{d}=[d_i]$ is the $N_n$-dimensional vector independent of $\mathbf{d}$. The boundary conditions are included in $\mathbf{d}$. The system of ODEs(49) is temporally integrated with the $\theta$-method (Knabner and Angermann, 2002) with $\theta > 0.5$. In this DFVM, the temporal integration of the advection-diffusion sub step is mathematically as

$$\mathbf{u}^{**} = (I-\Xi)^{-1}(\mathbf{u}^* + \mathbf{d})$$
(50)

where $I$ is the $N_n \times N_n$-dimensional identity matrix and the vectors $\mathbf{u}^{**}$ and $\mathbf{u}^*$ are given by $\mathbf{u}^* = [u_i^*]$ and $\mathbf{u}^{**} = [u_i^{**}]$, respectively. The inversion of the matrix $I-\Xi$ performed in Eq.(50) is mathematically justified because of the positivity conditions of the spatial discretization in the DFVM as proved in Yoshioka and Unami (2013). The recently developed point-implicit method (Kadioglu et al., 2015), which can be regarded as spatially local $\theta$-method, has been examined in preliminary investigation, but it computed severely smeared numerical solutions. The method is therefore not used in this article.

## 4. 2-D Dual-Finite Volume Method

The 1-D DFVM presented in the previous section can be straightforwardly extended to a 2-D counterpart, which is for simulating moisture dynamics in homogenous thin non-woven fibrous sheets where the moisture profiles are essentially 2-D.

### 4.1 Spatial Discretization
#### 4.1.1 Computational Mesh

The 2-D domain $\Omega$ is divided into triangular elements in the usual conforming manner (Ref). The elements and the nodes are indexed with the natural numbers. The total numbers of the elements and nodes are denoted by $N_e$ and $N_n$, respectively. The $i$ th node and the $k$ th element are denoted by $P_i$ and $\Omega_k$, respectively. The horizontal area of the element $\Omega_k$ is denoted by $|\Omega_k|$. The number of elements sharing the node $P_i$ is denoted by $v(i)$ and the $l$ th of them by $\Omega_{\kappa(i,j)}$. The element $\Omega_{\kappa(i,j)}$ has three nodes; one of them is the $i$ th node $P_i$ and the others are referred to as the $\mu(i,j,1)$ th node $P_{\mu(i,j,1)}$ and the $\mu(i,j,2)$ th node $P_{\mu(i,j,2)}$ in a counterclockwise manner. The number of nodes directly connected to the node $P_i$ is denoted by $\tau(i)$ and the $l$ th of them is denoted by the $\rho(i,l)$ th node $P_{\rho(i,l)}$. The length of the edge $P_i P_{\rho(i,l)}$ is denoted by $d_{i,l}$. A dual cell associated with the node $P_i$ is denoted by $S_i$, which is defined as the conventional 2-D Voronoi cell (Mishev, 1998) given by

$$S_i = \bigcup_{l=1}^{\tau(i)} S_{i,l} = \bigcup_{l=1}^{\tau(i)} \left\{ \mathbf{x} \middle| \text{dist}(\mathbf{x}, P_i) < \text{dist}(\mathbf{x}, P_{\rho(i,l)}) \right\} \tag{51}$$

where $\text{dist}(\cdot,\cdot)$ represents the conventional distance function between a couple of points located in the domain $\Omega$. The horizontal area of the dual cell $S_i$ is denoted by $|S_i|$. The length of the outer perimeter of the sub-cell $S_{i,l}$ of the dual cell $S_i$ is denoted by $L_{i,l}$. The unknown $u$ is attributed to the dual cells (or equivalently to the nodes).

### 4.2 Operator-Splitting Algorithm

The operator splitting algorithm used the 2-D DFVM is essentially similar to that in the 1-D counterpart. The 2-D DFVM solves Eq.(9) at each time step following Eq.(12). Discretization procedure of the evaporation and advection-diffusion two sub-steps at each time step is explained in the following subsections.

### 4.3 Evaporation Sub-step

At the evaporation sub-step, the equation to be solved at each node is given by

$$\frac{\partial u}{\partial t} = -E_s u^q, \tag{52}$$

which is essentially same with Eq.(13) except for that now the nodes are distributed in a 2-D domain. The temporal integration procedure in this sub-step also follows that of the 1-D counterpart.

### 4.4 Advection-diffusion Sub-step

At the advection-diffusion sub-step, the equation to be solved is given by

$$\frac{\partial u}{\partial t} + \frac{\partial F_1}{\partial x_1} + \frac{\partial F_2}{\partial x_2} = 0 \tag{53}$$

which can be formally regarded as a 2-D nonlinear advection-diffusion equation. A finite volume discretization is applied to Eq.(20) as presented in the following sub-sections.

#### 4.4.1   Finite Volume Formulation

Integrating the Eq.(53) in the dual cell $S_i$ with the application of the Gauss-Green theorem leads to the local integral relationship

$$\frac{\partial}{\partial t} \int_{S_i} u \, ds + \sum_{l=1}^{\tau(i)} L_{i,l} F_{i,l} = 0. \tag{54}$$

where $F_{i,l}$ is the numerical flux to be evaluated on the cell face between the dual cells $S_i$, $S_{\rho(i,l)}$. As in the 1-D DFVM, the three integrals appearing the first term in the left hand-side of Eq.(54) is evaluated as

$$\frac{\partial}{\partial t} \int_{S_i} u \, dx \approx |S_i| \frac{du_i}{dt}. \tag{55}$$

The numerical flux $F_{i,j}$ on the cell interface $\Gamma_{i,j}$ is evaluated with a fitting technique, which

uses the analytical solution $\tilde{u}$ of the two-point boundary value problem

$$\frac{\partial}{\partial s}\left(b_{\rho(i,l)}u - a_{\rho(i,l)}\frac{\partial u}{\partial s}\right) = 0 \tag{56}$$

along the edge $P_i P_{\rho(i,l)}$ with the local 1-D coordinate $s$ taken from the node $P_i$ to the node $P_{\rho(i,l)}$ and the nodal boundary conditions

$$u(x_i) = u_i, \quad u(x_{\rho(i,l)}) = u_{\rho(i,l)} \tag{57}$$

where the coefficients $a_{\rho(i,l)}$ and $b_{\rho(i,l)}$ are given by

$$a_{\rho(i,l)} = (m-p)\left[h_\varepsilon\left(\overline{u_{\rho(i,l)}}\right)\right]^{m-p-1} \tag{58}$$

and

$$b_{\rho(i,l)} = -\sin\alpha\left[h_\varepsilon\left(\overline{u_{\rho(i,l)}}\right)\right] n_{1,\rho(i,l)} \tag{59}$$

where $\overline{u_{\rho(i,l)}}$ represents the value of $u$ on the edge $P_i P_{\rho(i,l)}$ that is specified later and $n_{1,\rho(i,l)}$ is the $x_1$-component of the unit normal vector $(n_{1,\rho(i,l)}, n_{2,\rho(i,l)})$ having the same direction with the vector $\overrightarrow{P_i P_{\rho(i,l)}}$ connecting the two nodes $P_i$ and $P_{\rho(i,l)}$. The analytical solution $u$ is explicitly obtained as

$$u = \frac{\exp(p_{i,l})u_i - u_{\rho(i,l)}}{\exp(p_{i,l}) - 1} + \frac{u_{\rho(i,l)} - u_i}{\exp(p_{i,l}) - 1}\exp\left(\frac{n_{1,\rho(i,l)}b_{\rho(i,l)}}{a_{\rho(i,l)}}(s - s_i)\right), \tag{60}$$

which is used for directly evaluating the numerical flux $F_{i,l}$ as

$$\begin{aligned} F_{i,l} &= a_{\rho(i,l)}u - b_{\rho(i,l)}\frac{\partial u}{\partial s} \\ &= \frac{(m-p)n_{1,\rho(i,l)}}{d_{i,l}}\left[h_\varepsilon\left(\overline{u_{\rho(i,l)}}\right)\right]^{m-p-1} p_{i,l}\frac{\exp(p_{i,l})u_i - u_{\rho(i,l)}}{\exp(p_{i,l}) - 1} \end{aligned} \tag{61}$$

where the local Peclet number $p_{i,l}$, which in the present case depends on the solution $u$, is defined as

$$p_{i,l} = -\frac{n_{1,\rho(i,l)}d_{i,l}\sin\alpha}{m-p}h_\varepsilon\left(\overline{u_{\rho(i,l)}}\right)^p. \tag{62}$$

The numerical flux $F_{i,j}$ is fully specified if the value of $\overline{u_{\rho(i,l)}}$ is evaluated, which is analogously performed as in the 1-D DFVM, namely the CE, FU, and IS methods can be utilized. Assembling Eq.(54) with the above-mentioned spatial discretization procedure for every node yields a system of linear ODEs as in the 1-D case, which is performed with the $\theta$-method with $\theta > 0.5$.